\documentclass[12pt,a4paper,reqno]{amsart}
\usepackage{amsmath,amssymb,amsthm}
\usepackage{url}
\usepackage{fullpage}
\usepackage{comment}

\newcommand{\sameorder}{\asymp}
\newcommand{\dx}{\mathrm{d}}
\newcommand{\C}{\mathbb{C}}
\newcommand{\N}{\mathbb{N}}
\newcommand{\R}{\mathbb{R}}
\newcommand{\Y}{\mathbb{Y}}
\newcommand{\doublesum}{\mathop{\sum\sum}}
\newcommand{\Stilde}{\widetilde{S}}
\newcommand{\Odip}[2]{\mathcal{O}_{#1}\!\left(#2\right)\mathchoice{\!}{}{}{}}
\newcommand{\Odi}[1]{\Odip{}{#1}}
\newcommand{\odip}[2]{{o}_{#1}\!\left(#2\right)\mathchoice{\!}{}{}{}}
\newcommand{\odi}[1]{\odip{}{#1}}

\allowdisplaybreaks

\renewcommand{\qedsymbol}{$\square$}
\newenvironment{Proof}[1][Proof]{\par\noindent\textbf{#1.}~}
{\hfill\qedsymbol\smallskip\par}
\newtheorem{Theorem}{Theorem}
\newtheorem{Lemma}{Lemma}

\begin{document}

\title{A Ces\`aro Average of Goldbach numbers}
\date{}
\author{Alessandro Languasco \& Alessandro Zaccagnini}
%
\subjclass[2010]{Primary 11P32; Secondary 44A10}
\keywords{Goldbach-type theorems, Laplace transforms, Ces\`aro averages}
\begin{abstract}
Let $\Lambda$ be the von Mangoldt function and
\(
r_G(n) =   \sum_{m_1 + m_2 = n} \Lambda(m_1) \Lambda(m_2)
\)
be the counting function for the Goldbach numbers. Let $N \geq 2$
be an integer. We prove that
\begin{align*}
  &\sum_{n \le N} r_G(n) \frac{(1 - n/N)^k}{\Gamma(k + 1)}
  =
  \frac{N^{2}}{\Gamma(k + 3)}
  -
  2
  \sum_{\rho} \frac{\Gamma(\rho)}{\Gamma(\rho + k + 2)} N^{\rho+1} \\
  &\qquad+
  \sum_{\rho_1} \sum_{\rho_2}
    \frac{\Gamma(\rho_1) \Gamma(\rho_2)}{\Gamma(\rho_1 + \rho_2 + k + 1)}
    N^{\rho_1 + \rho_2}
  +
  \Odip{k}{N^{1/2}},
\end{align*}
for $k > 1$, where $\rho$, with or without subscripts, runs over
the non-trivial zeros of the Riemann zeta-function $\zeta(s)$.
\end{abstract}
\maketitle

\section{Introduction}

We continue our recent work on the number of representations of an
integer as a sum of primes.
In \cite{LanguascoZ2012a} we studied the \emph{average} number of
representations of an integer as a sum of two primes, whereas in
\cite{LanguascoZ2012b} we considered individual integers.
In this paper we study a Ces\`aro weighted \emph{explicit} formula for
Goldbach numbers and the goal is similar to the one in
\cite{LanguascoZ2012a}, that is, we want to obtain the expected main
term and one or more terms that depend explicitly on the zeros of the
Riemann zeta-function, with a small error.
Letting
\[
  r_G(n)
  =
  \sum_{m_1 + m_2 = n} \Lambda(m_1) \Lambda(m_2),
\]
the main result of the paper is the following theorem.

\begin{Theorem}
\label{CesaroG-average}
Let $N$ be a positive integer.
We have
\begin{align}
\notag
  &\sum_{n \le N} r_G(n) \frac{(1 - n/N)^k}{\Gamma(k + 1)}
  =
  \frac{N^{2}}{\Gamma(k + 3)}
  -
  2
  \sum_{\rho} \frac{\Gamma(\rho)}{\Gamma(\rho + k + 2)} N^{\rho+1} \\
  &\qquad+
\label{expl-form-Goldbach}
  \sum_{\rho_1} \sum_{\rho_2}
    \frac{\Gamma(\rho_1) \Gamma(\rho_2)}{\Gamma(\rho_1 + \rho_2 + k + 1)}
    N^{\rho_1 + \rho_2}
  +
  \Odip{k}{N^{1/2}},
\end{align}
for $k > 1$, where $\rho$, with or without subscripts, runs over
the non-trivial zeros of the Riemann zeta-function $\zeta(s)$.
\end{Theorem}

We remark that the double series over zeros in \eqref{expl-form-Goldbach}
converges absolutely for $k > 1 / 2$, and it seems reasonable to
believe that the stated equality holds for the same values of $k$,
possibly with a weaker error term, although the bound $k>1$ appears in
several places of the proof and it seems to be the limit of the
method.

The result in \cite{LanguascoZ2012a} is the case $k = 0$ of
\eqref{expl-form-Goldbach} under the assumption of the Riemann
Hypothesis (RH); there we only get the first sum over zeros
and the error term is $\Odi{N (\log N)^3}$. The proof in
\cite{LanguascoZ2012a} depends on RH in just one place; it is 
not hard to get an unconditional version of such a result with 
an error term  $\odi{N^{2}}$.
The technique here is completely different and RH has no consequences
on the lower bound for the size of $k$.

Similar averages of arithmetical functions are common in the 
literature, see, e.g., Chan\-dra\-sekharan-Narasimhan
\cite{ChandrasekharanN1961} and Berndt \cite{Berndt1975} 
who built on earlier classical works.
In their setting the generalized Dirichlet series associated to the
arithmetical function satisfies a suitable functional equation and
this leads to an asymptotic formula containing Bessel functions of
real order.
In our case we have no functional equation, and Bessel functions are
naturally replaced by Gamma functions; in fact we plan to develop
further the present technique to deal with the cases
$p_1^{\ell_{1}}+p_2^{\ell_{2}}$ and $p+m^{2}$, where Bessel functions
with complex order arise; we expect many technical complications.

The most interesting explicit formula in Goldbach's problem has been
recently given by Pintz \cite{Pintz2006}.
It is too complicated to be reproduced here, but we remark that in his
formula, which deals with individual values of $r_G(n)$, the summation
is over zeros of suitable Dirichlet $L$-functions, whereas in an
average problem like the present one, only the zeros of the Riemann
$\zeta$-function are relevant.
The same phenomenon occurs in our papers \cite{LanguascoZ2012a} and
\cite{LanguascoZ2012b}.

The method we will use is based on a formula due to Laplace
\cite{Laplace1812}, namely
\begin{equation}
\label{Laplace-transf}
  \frac 1{2 \pi i}
  \int_{(a)} v^{-s} e^v \, \dx v
  =
  \frac1{\Gamma(s)},
\end{equation}
where $\Re(s) > 0$ and $a > 0$, see, e.g., formula 5.4(1) on page 238
of \cite{ErdelyiMOT1954a}.
In the following we will need the general case of \eqref{Laplace-transf}
which can be found in de Azevedo Pribitkin \cite{Azevedo2002}, 
formulae (8) and (9):
\begin{equation}
\label{Laplace-eq-1}
  \frac1{2 \pi}
  \int_{\R} \frac{e^{i D u}}{(a + i u)^s} \, \dx u
  =
  \begin{cases}
    \dfrac{D^{s - 1} e^{- a D}}{\Gamma(s)}
    & \text{if $D > 0$,} \\
    0
    & \text{if $D < 0$,}
  \end{cases}
\end{equation}
which is valid for $\sigma = \Re(s) > 0$ and $a \in \C$ with
$\Re(a) > 0$, and
\begin{equation}
\label{Laplace-eq-2}
  \frac1{2 \pi}
  \int_{\R} \frac 1{(a + i u)^s} \, \dx u
  =
  \begin{cases}
    0     & \text{if $\Re(s) > 1$,} \\
    1 / 2 & \text{if $s = 1$,}
  \end{cases}
\end{equation}
for $a \in \C$ with $\Re(a) > 0$.
Formulae \eqref{Laplace-eq-1}-\eqref{Laplace-eq-2} enable us to write
averages of arithmetical functions by means of line integrals as we
will see in \S\ref{settings} below.
We recall that Walfisz, see \cite[Ch.~X]{Walfisz1957}, replaced
\eqref{Laplace-eq-1}-\eqref{Laplace-eq-2} with the following
particular case
\[
  \frac1{2 \pi i}
  \int_{(a)} e^{x \omega} \frac{\dx \omega}{\omega^{\ell + 1}}
  =
  \begin{cases}
    x^\ell / \ell! & \text{if $x > 0$,} \\
    0              & \text{if $x \le 0$,}
  \end{cases}
\]
which is valid for $\ell \in \N$ with $\ell \ge 1$, $a>0$,
and $x \in \R$.

We combine this approach with line integrals with the classical
methods dealing with infinite sums, exploited by Hardy \& Littlewood
(see \cite{HardyL1916} and \cite{HardyL1923}) and by Linnik
\cite{Linnik1946}.
In particular, in \S2.5 of \cite{HardyL1916} there is a sort of
``explicit formula'' for a function related to $\psi(x) - x$.

We thank A.~Perelli and J.~Pintz for several conversations on this topic.

\section{Settings}
\label{settings}

Let
\begin{equation}
\label{Stilde-def}
  \Stilde(z)
  =
  \sum_{m \ge 1} \Lambda(m) e^{- m z},
\end{equation}
where $z = a + i y$ with $y \in \R$ and real $a > 0$.
We recall that the Prime Number Theorem (PNT) is equivalent, via
Lemma~\ref{Linnik-lemma2} below, to the statement
\begin{equation}
\label{PNT-equiv}
  \Stilde(a)
  \sim
  a^{-1}
  \qquad\text{for $a \to 0+$,}
\end{equation}
which is classical: for the proof see for instance Lemma~9 in Hardy \&
Littlewood \cite{HardyL1923}.
By \eqref{Stilde-def} we have
\[
  \Stilde(z)^2
  =
  \sum_{n \ge 1} r_G(n) e^{-n z}.
\]
Hence, for $N \in \N$ with $N > 0$ and $a > 0$ we have
\begin{equation}
\label{starting-point}
 \frac 1{2 \pi i}
  \int_{(a)} e^{N z} z^{- k - 1} \Stilde(z)^2 \, \dx z
  =
  \frac 1{2 \pi i}
  \int_{(a)} e^{N z} z^{- k - 1}
    \sum_{n \ge 1} r_{G}(n) e^{- n z} \, \dx z.
\end{equation}
Since
\[
  \sum_{n \ge 1} \vert r_G(n) e^{- n z} \vert 
  =
  \Stilde(a)^2
  \sameorder
  a^{-2}
\]
by \eqref{PNT-equiv}, where $f\sameorder g$ means $g \ll f \ll g$, we
can exchange the series and the line integral in
\eqref{starting-point} provided that $k>0$.
In fact, if $z = a + i y$, taking into account the estimate
\begin{equation}
\label{z^-1}
  \vert z \vert^{-1}
  \sameorder
  \begin{cases}
    a^{-1}   &\text{if $\vert y \vert \le a$,} \\
    \vert y \vert^{-1} &\text{if $\vert y \vert \ge a$,}
  \end{cases}
\end{equation}
we have
\[
  \vert e^{N z} z^{- k - 1}\vert 
  \sameorder
   e^{N a}
  \begin{cases}
    a^{- k - 1} &\text{if $\vert y \vert \le a$,} \\
    \vert y \vert^{- k - 1} &\text{if $\vert y \vert \ge a$,}
  \end{cases}
\]
and hence, recalling \eqref{PNT-equiv}, we obtain
\begin{align*}
  \int_{(a)} \vert e^{N z} z^{- k - 1} \vert \,
    \Bigl\vert
      \sum_{n \ge 1} r_{G}(n) e^{- n z}
    \Bigr\vert \, \vert \dx z \vert
  &\ll
  a^{-2} e^{N a}
  \Bigl[
    \int_{-a}^a a^{- k - 1} \, \dx y
    +
    2
    \int_a^{+\infty} y^{- k - 1} \, \dx y
  \Bigl] \\
  &=
  2 a^{-2} e^{N a}
  \Bigl( a^{-k} + \frac{a^{-k}}k \Bigr),
\end{align*}
but only for $k > 0$.
Using \eqref{Laplace-eq-1} for $n \ne N$ and \eqref{Laplace-eq-2} for
$n = N$, we see that the right-hand side of \eqref{starting-point} is
\begin{align*}
  &=
  \sum_{n \ge 1} r_{G}(n)
    \Bigl[
      \frac 1{2 \pi i}
      \int_{(a)} e^{(N - n) z} z^{- k - 1} \, \dx z
    \Bigr] 
  =
  \sum_{n \le N} r_{G}(n) \frac{(N - n)^k}{\Gamma(k + 1)}
\end{align*}
for $k > 0$.

\paragraph{\textbf{Remark.}}
It is important to notice that the previous computation reveals
that we can not get rid of the Ces\`aro weight in our method since,
for $k = 0$, it is not clear if the integral at the right hand side
of \eqref{starting-point} converges absolutely or not.
In fact, if we could prove \eqref{expl-form-Goldbach} for $k = 0$,
assuming the RH we could easily derive the main result of
\cite{LanguascoZ2012a} with an error term $\Odi{N}$, and this seems to
be quite unreachable in the present state of knowledge.
See the concluding remarks in the latter paper for an explanation.

\medskip
Summing up
\begin{equation*}
  \sum_{n \le N}
    r_{G}(n) \frac{(N - n)^k}{\Gamma(k + 1)}
  =
 \frac 1{2 \pi i}
  \int_{(a)} e^{N z} z^{- k - 1} \Stilde(z)^2 \, \dx z,
\end{equation*}
where $N \in \N$ with $N > 0$, $a>0$ and $k > 0$.
This is the fundamental relation for the method.

\section{Inserting zeros}

In this section we need $k > 1$.
By Lemma \ref{Linnik-lemma2} below we have
\[
  \Stilde(z)
  =
  \frac 1z
  -
  \sum_{\rho} z^{-\rho} \Gamma(\rho)
  +
  E(y,a)
\]
where $E(y, a)$ satisfies \eqref{expl-form-err-term-strong}.
Hence
\[
  \Stilde(z)^2
  =
  \Bigl( \frac 1z - \sum_{\rho} z^{-\rho} \Gamma(\rho) \Bigr)^2
  +
  E(y,a)^2
  +
  2 E(y,a)
  \Bigl( \frac 1z - \sum_{\rho} z^{-\rho} \Gamma(\rho) \Bigr).
\]
We have
\[
  \Bigl\vert \frac 1z - \sum_{\rho} z^{-\rho} \Gamma(\rho) \Bigr\vert
  =
  \bigl\vert \Stilde(z) - E(y,a) \bigr\vert
  \le
  \Stilde(a) + \bigl\vert E(y,a) \bigr\vert
  \ll
  a^{-1}
  +
  \bigl\vert E(y,a) \bigr\vert
\]
by \eqref{PNT-equiv} again, so that
\begin{equation}
\label{squaring-out}
  \Stilde(z)^2
  =
  \Bigl( \frac 1z - \sum_{\rho} z^{-\rho} \Gamma(\rho) \Bigr)^2
  +
  \Odi{\bigl\vert  E(y,a) \bigr\vert^2
  +
  \bigl\vert  E(y,a) \bigr\vert a^{-1}}.
\end{equation}
Recalling \eqref{z^-1} and \eqref{expl-form-err-term-strong}, we have
\begin{align*}
  \int_{(a)} \bigl \vert E(y,a) \bigr\vert^2
    \vert e^{N z} \vert \, \vert z \vert^{- k - 1} \, \vert \dx z \vert
  &\ll
  e^{N a}
  \int_0^a a^{- k} \, \dx y
  +
  e^{N a}
  \int_a^{+\infty} y^{- k}(1+ \log^2(y/a))^2 \, \dx y \\
  &\ll_k
  e^{N a} a^{- k + 1} +  e^{N a} a^{- k + 1}
  \int_1^{+\infty} v^{- k} (1 + \log^2v)^2 \, \dx v \\
  &\ll_k
  e^{N a} a^{- k + 1}.
\end{align*}
Choosing $a = 1 / N$, the error term is $\ll_k N^{k - 1}$ for $k > 1$.
For $a = 1 / N$, by \eqref{z^-1} and \eqref{expl-form-err-term-strong},
the second remainder term in \eqref{squaring-out} is
\begin{align*}
  &\ll
  N
  \int_{(1/N)} \vert E(y,1/N) \vert\vert e^{N z}\vert
    \vert z\vert ^{- k - 1} \, \vert \dx z \vert \\
  &\ll
  N
  \int_0^{1/N} N^{ k + 1/2} \, \dx y
  +
  N
  \int_{1/N}^{+\infty} y^{- k - 1/2} \log^2(Ny) \, \dx y \\
  &\ll
  N^{k + 1/2} 
  +
  N^{k + 1/2} 
  \int_{1}^{+\infty}  v^{- k - 1/2} \log^2 v \, \dx v 
  \ll_k
  N^{k + 1/2}.
\end{align*}
With a little effort we can give an explicit dependence on $k$ for the
implicit constants in the last two estimates, showing that the
condition $k>1$ is indeed necessary.

Hence, by \eqref{starting-point} we have
\begin{align}
\notag
  \sum_{n \le N} r_G(n) &\frac{(N - n)^k}{\Gamma(k + 1)}
  =
  \frac 1{2 \pi i}
  \int_{(1/N)}
    e^{N z} z^{- k - 1}
    \Bigl( \frac 1z - \sum_{\rho} z^{-\rho} \Gamma(\rho) \Bigr)^2
    \, \dx z
  +
  \Odip{k}{N^{k + 1/2}} \\
  \notag
  &=
  \frac 1{2 \pi i}
  \int_{(1/N)} e^{N z} z^{- k - 3} \, \dx z
  -
  \frac 1{\pi i}
  \int_{(1/N)} e^{N z} z^{- k - 2}
    \sum_{\rho} z^{-\rho} \Gamma(\rho) \, \dx z \\
  \label{inserting-zeros}
  &\quad+
  \frac 1{2 \pi i}
  \int_{(1/N)}
    e^{N z} z^{- k - 1}
    \sum_{\rho_1} \sum_{\rho_2}
      z^{- \rho_1 - \rho_2} \Gamma(\rho_1) \Gamma(\rho_2)
    \, \dx z
  +
  \Odip{k}{N^{k + 1/2}}.
\end{align}
Interchanging the series with the integrals 
(see \S\ref{first-exchange}-\ref{exchange-double-sum-rhos} for a proof
that this is permitted when $k>1$), by \eqref{inserting-zeros} we get
that
\begin{align*}
  \sum_{n \le N} r_G(n) \frac{(N - n)^k}{\Gamma(k + 1)}
  &=
  \frac 1{2 \pi i}
  \int_{(1/N)} e^{N z} z^{- k - 3} \, \dx z
  -
  \frac 1{\pi i}
  \sum_{\rho} \Gamma(\rho)
  \int_{(1/N)} e^{N z} z^{- k - 2 - \rho} \, \dx z \\
  &\qquad+
  \frac 1{2 \pi i}
  \sum_{\rho_1} \sum_{\rho_2} \Gamma(\rho_1) \Gamma(\rho_2)
  \int_{(1/N)} e^{N z} z^{- k - 1 - \rho_1 - \rho_2} \, \dx z
  +
  \Odip{k}{N^{k + 1/2}} \\
  &=
  I_1 + I_2 + I_3 + \Odip{k}{N^{k + 1/2}},
\end{align*}
say.

\subsection{Evaluation of $I_1$}

Using \eqref{Laplace-transf} and putting $s = N z$, we immediately get
\[
  I_1
  =
  \frac{N^{k + 2}}{2 \pi i}
  \int_{(1)} e^s s^{- k - 3} \, \dx s
  =
  \frac{N^{k + 2}}{\Gamma(k + 3)}.
\]

\subsection{Evaluation of $I_2$}

Putting $s = N z$ and by \eqref{Laplace-transf} again, we have
\begin{align*}
  I_2
  &=
  -
  \frac 1{\pi i}
  \sum_{\rho} \Gamma(\rho) N^{k + \rho + 1}
  \int_{(1)} e^{s} s^{- k - 2 - \rho} \, \dx s
  =
  - 2
  \sum_{\rho} \frac{\Gamma(\rho)}{\Gamma(\rho + k + 2)} N^{k + \rho + 1}. 
\end{align*}

\subsection{Evaluation of $I_3$}

As above, using \eqref{Laplace-transf} and putting $s = N z$, we get
\begin{align*}
  I_3
  &=
  \frac 1{2 \pi i}
  \sum_{\rho_1} \sum_{\rho_2} \Gamma(\rho_1) \Gamma(\rho_2) N^{k + \rho_1 + \rho_2}
  \int_{(1)} e^s s^{- k - 1 - \rho_1 - \rho_2} \, \dx s \\
  &=
  \sum_{\rho_1} \sum_{\rho_2}
    \frac{\Gamma(\rho_1) \Gamma(\rho_2)}{\Gamma(\rho_1 + \rho_2 + k + 1)}
    N^{k + \rho_1 + \rho_2}.
\end{align*}
Combining the previous relations, we finally get
\begin{align}
\notag
  &\sum_{n \le N} r_G(n) \frac{(N - n)^k}{\Gamma(k + 1)}
  =
  \frac{N^{k + 2}}{\Gamma(k + 3)}
  -
  2 N^{k + 1}
  \sum_{\rho} \frac{\Gamma(\rho)}{\Gamma(\rho + k + 2)} N^{\rho} \\
  &\qquad+
\label{expl-form-Goldbach-bis}
  N^{k} \sum_{\rho_1} \sum_{\rho_2}
    \frac{\Gamma(\rho_1) \Gamma(\rho_2)}{\Gamma(\rho_1 + \rho_2 + k + 1)}
    N^{\rho_1 + \rho_2} 
  +
  \Odip{k}{N^{k + 1/2}}
\end{align}
for $k > 1$.
The proof that the double sum over zeros converges absolutely for
$k > 1 / 2$ is given in \S\ref{sec:double-sum} below.
Theorem \ref{CesaroG-average} follows dividing
\eqref{expl-form-Goldbach-bis} by $N^{k}$.

\section{Lemmas}
We recall some basic facts in complex analysis.
First, if $z = a + i y$ with $a > 0$, we see that for complex $w$ we
have
\begin{align*}
  z^{-w}
  &=
  \vert z \vert^{-w} \exp( - i w \arctan(y / a)) \\
  &=
  \vert z \vert^{-\Re(w) - i \Im(w)} \exp( (- i \Re(w) + \Im(w)) \arctan(y / a))
\end{align*}
so that
\begin{equation}
\label{z^w}
  \vert z^{-w} \vert
  =
  \vert z \vert^{-\Re(w)} \exp(\Im(w) \arctan(y / a)).
\end{equation}
We also recall that, uniformly for $x \in [x_1, x_2]$, with $x_1$ and
$x_2$ fixed, and for $|y| \to +\infty$, by the Stirling formula we have
\begin{equation}
\label{Stirling}
  \vert \Gamma(x + i y) \vert
  \sim
  \sqrt{2 \pi}
  e^{- \pi |y| / 2} |y|^{x - 1 / 2},
\end{equation}
see, e.g., Titchmarsh \cite[\S4.42]{Titchmarsh1988}.

We will need the Hardy-Littlewood-Linnik formula (see, e.g., Languasco
\& Zaccagnini \cite{LanguascoZ2012b}): we
notice that here $y \in \R$, while in \cite{LanguascoZ2012b} we had
the restricted range $y \in [-1/2, 1/2]$.
Hence there are some modifications to be made.
We will follow the proof in Linnik \cite{Linnik1946}
(see also eq.~(4.1) of \cite{Linnik1952}).

\begin{Lemma} 
\label{Linnik-lemma2}
Let $z = a + iy$, where $a > 0$ and $y \in \R$.
Then
\begin{equation*}
  \widetilde{S}(z)
  =
  \frac{1}{z}
  -
  \sum_{\rho}z^{-\rho} \Gamma(\rho)
  +
  E(a,y)
\end{equation*}
where $\rho = \beta + i\gamma$ runs over the non-trivial zeros of
$\zeta(s)$ and
\begin{equation}
\label{expl-form-err-term-strong}
  E(a,y)
  \ll
  \vert z \vert^{1/2}
  \begin{cases}
    1 & \text{if $\vert y \vert \leq a$} \\
    1 +\log^2 (\vert y\vert/a) & \text{if $\vert y \vert > a$.}
  \end{cases}
\end{equation}
\end{Lemma}

\begin{Proof}
Following the line of Hardy and Littlewood, see
\cite[\S2.2]{HardyL1916}, \cite[Lemma 4]{HardyL1923} and of \S4 in
Linnik \cite{Linnik1946}, we have that
\begin{equation}
\label{Mellin2}
  \widetilde{S}(z)
  =
  \frac{1}{z}
  -
  \sum_{\rho}z^{-\rho}\Gamma(\rho)
  -
  \frac{\zeta'}{\zeta}(0)
  -
  \frac{1}{2\pi i}
  \int_{(-1/2)} 
  \frac{\zeta'}{\zeta}(w) \Gamma(w)z^{-w} \, \dx w.
\end{equation}
Now we estimate the integral in \eqref{Mellin2}. 
Let $c>0$ be a positive constant to be chosen later.
Writing $w=-1/2+it$, we have
$\vert (\zeta'/\zeta)(w)\vert \ll \log (\vert t \vert +2)$,
$\vert z^{-w} \vert = \vert z \vert^{1 / 2} \exp(t \arctan( y / a ))$
by \eqref{z^w} and, for $\vert t\vert > c$,
$\Gamma(w) \ll \vert t \vert^{-1} \exp(-\frac{\pi}{2}\vert t \vert)$
by \eqref{Stirling}.
Letting $L_{c}=\{-1/2+it: \vert t \vert > c\}$  we have
\[
\int_{L_{c}} 
\frac{\zeta'}{\zeta}(w)  \Gamma(w) z^{-w} \, \dx w
  \ll
   \vert z \vert^{1/2}
   \int_{L_{c}} 
     \frac{\log \vert t \vert}{\vert t \vert}
     \exp
     \Bigl(
     -\frac{\pi}{2}\vert t \vert + t \arctan(  y / a )
     \Bigr)
     \, \dx t.
\]
If $t y\leq 0$ we call $\eta$ the quantity 
$\frac{\pi}{2}+ \vert\arctan(  y / a )\vert \in [\pi/2, \pi)$.
If $\vert y \vert \leq a$ we define $\eta$ as 
$\frac{\pi}{2} - \arctan( y / a )>\frac{\pi}{2} - \arctan(1)=
\frac{\pi}{4}$.
In the remaining case ($\vert y \vert > a$ and $ty > 0$)
we set $\eta=\arctan(a/ \vert y\vert) \gg a/\vert y \vert$.
Now fix $c$ such that $c\eta<1$ (e.g., $c=1/\pi$ is allowed). 
Letting $u = \eta t$, we get
\begin{align}
  \notag
  \int_{L_{c}}
  \frac{\zeta'}{\zeta}(w)  \Gamma(w) z^{-w} \, \dx w
  &\ll
  \vert z \vert^{1/2}
  \int_c^{+\infty} e^{-\eta t} \, \frac{\log t}t \, \dx t
  =
  \vert z \vert^{1/2}
  \int_{c \eta}^{+\infty} e^{-u} \, \frac{\log(u / \eta)}u \, \dx u
  \\
  \notag
  &=
  \vert z \vert^{1/2}
  \int_{c \eta}^{+\infty} e^{-u} \, \frac{\log u}u \, \dx u
  +
  \vert z \vert^{1/2}
  \log(1 / \eta)
  \int_{c \eta}^{+\infty} e^{-u} \frac{\dx u}u \\
\label{splitting}
  &=
  J_1 + J_2.
\end{align}
We remark that $0 \le u^{-1} \log u \le e^{-1}$ for $u \ge 1$, since
the maximum of $u^{-1} \log u$ is attained at $u = e$.
Since
\[
  0
  \le
  \int_1^{+\infty} e^{-u} \, \frac{\log u}u \, \dx u
  \le
  e^{-1}
  \int_1^{+\infty} e^{-u} \, \dx u
  \ll
  1
\]
and
\[
  \Bigl\vert
    \int_{c \eta}^1 e^{-u} \, \frac{\log u}u \, \dx u
  \Bigr\vert
  \le
  \int_{c \eta}^1 \frac{-\log u}u \, \dx u
  =
  \Bigl[ - \frac12 \log^2 u
  \Bigr]_{c \eta}^1
  \ll
  \log^2 (1 / \eta)
\]
we have that $J_{1} \ll \vert z \vert^{1/2} \log^2 (1 / \eta)$.
For $J_{2}$ it is sufficient to remark that
\[
  0
  \le
  J_2
  \le
  \vert z \vert^{1/2}
  \log(1 / \eta)
  \Bigl(
    \int_{c \eta}^1 \frac{\dx u}u
    +
    \int_1^{+\infty} e^{-u} \dx u
  \Bigr)
  \ll
  \vert z \vert^{1/2}
  \log^2 (1 / \eta).
\]
Inserting the last two estimates in \eqref{splitting}, recalling the
definition of $\eta$ and remarking that the integration over
$\vert t \vert \leq c$ gives immediately a contribution $\ll 1$, we
obtain that the integral in \eqref{Mellin2} is dominated by the right
hand side of \eqref{expl-form-err-term-strong} and the lemma is
proved.
\end{Proof}
In the next sections we will need to perform several times a set of 
similar computations; so we collected them in the following two lemmas.
\begin{Lemma}
\label{series-int-zeros}
Let $\beta + i \gamma$ run over the non-trivial zeros of the Riemann
zeta-function and $\alpha > 1$ be a parameter.
The series
\[
  \sum_{\rho \colon \gamma > 0}
  \gamma^{\beta-1/2}
    \int_1^{+\infty} \exp\Bigl( - \gamma \arctan\frac 1u \Bigr)
      \frac{\dx u}{u^{\alpha+\beta}}
\]
converges provided that $\alpha > 3/2$.
For $\alpha \le 3/2$ the series does not converge.
The result remains true if we insert in the integral a factor
$(\log u)^c$, for any fixed $c \ge 0$.
\end{Lemma}

\begin{Proof}
Setting $y = \arctan(1 / u)$, for any real $\gamma > 0$ we have
\begin{align*}
  \int_1^{+\infty} \exp\Bigl( -\gamma \arctan\frac 1u \Bigr)
    \frac{\dx u}{u^{\alpha+\beta}}
  &=
  \int_0^{\pi / 4}
    \exp(-\gamma y) \,
    \frac{(\sin y)^{\alpha+\beta - 2}}{(\cos y)^{\alpha+\beta}} \, \dx y\\
  &\ll_\alpha
  \int_0^{\pi / 4}
    \exp(-\gamma y) \, y^{\alpha+\beta - 2} \, \dx y \\
  &=
  \gamma^{1 - \alpha - \beta}
  \int_0^{\pi \gamma / 4}
    \exp(-w) \, w^{\alpha+\beta - 2} \, \dx w \\
  &\le
  \gamma^{1 - \alpha - \beta} \,
  \Bigl(\Gamma(\alpha-1)+\Gamma(\alpha)\Bigr),
\end{align*}
since $0 < \beta < 1$.
This shows that the series over $\gamma$ converges for $\alpha > 3/2$.
For $\alpha = 3/2$ essentially the same computation shows that the
integral is $\gg \gamma^{-1/2 - \beta}$ and it is well known that in
this case the series over zeros diverges.
The other assertions are proved in the same way.
\end{Proof}

\begin{Lemma}
\label{series-int-zeros-alt-sign}
Let $\alpha > 1$, $z=a+iy$, $a\in(0,1)$ and $y\in \R$.
Let further $\rho=\beta+i\gamma$ run over the non-trivial zeros of
the Riemann zeta-function.  We have
\[
  \sum_{\rho}
    \vert \gamma\vert ^{\beta-1/2}
    \int_{\Y_1 \cup \Y_2} \exp\Bigl(\gamma \arctan\frac{y}{a} - \frac\pi2 \vert \gamma \vert\Bigr)
      \frac{\dx y}{\vert z \vert ^{\alpha+\beta}}
  \ll_{\alpha}
  a^{-\alpha},
\]
where $\Y_1=\{y\in \R\colon y\gamma \leq 0\}$ and 
$\Y_2=\{y\in [-a,a] \colon y\gamma > 0\}$.
The result remains true if we insert in the integral a factor
$(\log (\vert y\vert /a))^c$, for any fixed $c \ge 0$.
\end{Lemma}

\begin{Proof}
We first work on $\Y_1$.
By symmetry, we may assume that $\gamma > 0$.
For $y \in(-\infty, 0]$ we have
$\gamma \arctan(y/a) -\frac \pi2 \vert \gamma\vert  \le - \frac \pi2 \vert \gamma\vert $
and hence the quantity we are estimating becomes
\[
  \sum_{\rho \colon \gamma > 0} 
  \gamma^{\beta-1/2}
  \exp\Bigl( -\frac \pi2 \gamma \Bigr)
    \int_{-\infty}^0 \frac{\dx y}{\vert z \vert ^{\alpha+\beta}}
  \ll_\alpha
  \sum_{\rho \colon \gamma > 0} 
  \gamma^{\beta-1/2}
  \exp\Bigl( -\frac \pi2 \gamma \Bigr)
  a^{1-\alpha-\beta}
  \ll_{\alpha}
  a^{-\alpha},
\]
using $0<\beta<1$, standard zero-density estimates and \eqref{z^-1}.
We consider now the integral over $\Y_2$.
Again by symmetry we can assume that $\gamma > 0$ and so we get
\begin{align*}
  \sum_{\rho \colon \gamma > 0}
   \gamma ^{\beta-1/2}
    \int_0^{a} \exp\Bigl( \gamma (\arctan \frac{y}{a} -\frac \pi2) \Bigr)
     \frac{\dx y}{\vert z \vert ^{\alpha+\beta}}
  &\ll
  \sum_{\rho \colon \gamma > 0}
  \gamma ^{\beta-1/2}
  \exp\Bigl( -\frac \pi4 \gamma \Bigr)
  \int_0^{a} \frac{\dx y}{\vert z \vert ^{\alpha+\beta}} \\
  &\ll_\alpha
  \sum_{\rho \colon \gamma > 0}
  \gamma ^{\beta-1/2} \exp\Bigl(  -\frac \pi4 \gamma \Bigr)
  a^{1-\alpha-\beta} 
  \ll_\alpha
  a^{-\alpha}
\end{align*}
arguing as above.
The other assertions are proved in the same way.
\end{Proof}

\section{Interchange of the series over zeros with the line integral}
\label{first-exchange}

We need $k>1/2$ in this section. We need to establish the convergence of
\begin{equation}
\label{conv-integral}
  \sum_{\rho}
    \vert \Gamma(\rho) \vert
    \int_{(1/N)} \vert e^{N z} \vert \, \vert z \vert ^{- k - 1} \,
      \vert z^{- \rho} \vert \, \vert \dx z \vert.
\end{equation}
By \eqref{z^w} and the Stirling formula \eqref{Stirling}, we are left
with estimating
\begin{equation}
\label{conv-integral-1}
  \sum_{\rho}
    \vert \gamma\vert ^{\beta - 1/2}
    \int_{\R} \exp\Bigl( \gamma \arctan(N y) -\frac \pi2 \vert \gamma\vert \Bigr)
             \frac{\dx y}{\vert z \vert ^{k + 1 +\beta}}. 
\end{equation}
We have just to consider the case $\gamma y >0$, $\vert y \vert > 1/N$
since in the other cases the total contribution is  $\ll_k N^{k + 1 }$
by Lemma \ref{series-int-zeros-alt-sign} with $\alpha=k+1$ and $a=1/N$.
By symmetry, we may assume that $\gamma > 0$.
We have that the integral in \eqref{conv-integral-1} is
\begin{align*}
  &\ll
  \sum_{\rho \colon \gamma > 0}
    \gamma ^{\beta - 1/2}
    \int_{1 / N}^{+\infty} \exp\Bigl( - \gamma \arctan\frac 1{N y} \Bigr)
      \frac{\dx y}{y^{k + 1 +\beta}} \\
  &=
  N^{k}
  \sum_{\rho \colon \gamma > 0}
    N^{\beta}
    \gamma ^{\beta - 1/2}
    \int_1^{+\infty} \exp\Bigl( - \gamma \arctan\frac 1u \Bigr)
      \frac{\dx u}{u^{k + 1 +\beta}}.
\end{align*}
For $k > 1 / 2$ this is $\ll_k N^{k + 1}$ by
Lemma~\ref{series-int-zeros}. This implies
that the integrals in \eqref{conv-integral-1} and in \eqref{conv-integral} 
are both  $\ll_k N^{k + 1}$ and
hence this exchange step is fully justified.

\section{Interchange of the double series over zeros with the line integral} 
\label{exchange-double-sum-rhos}

We need $k>1$ in this section.
Arguing as in \S \ref{first-exchange}, we first need to establish the
convergence of
\begin{equation}
\label{conv-integral-2}
  \sum_{\rho_{1}}
    \vert \Gamma(\rho_{1})\vert  
    \int_{(1/N)}  
      \vert \sum_{\rho_{2}} \Gamma(\rho_{2})z^{- \rho_{2}} \vert
      \vert e^{N z} \vert \,
      \vert z \vert ^{- k - 1} \, \vert z^{- \rho_{1}} \vert \, \vert \dx z \vert.
\end{equation}
Using the PNT and \eqref{expl-form-err-term-strong}, we first remark
that
\begin{align}
\notag
  \Bigl\vert \sum_{\rho} z^{-\rho} \Gamma(\rho) \Bigr\vert 
  &=
  \bigl\vert  \Stilde(z) - \frac 1z -E(y,\frac{1}{N}) 
  \bigr\vert 
  \ll
  N  + \frac{1}{\vert z\vert}
  +
  \bigl\vert   E(y,\frac{1}{N}) \bigr\vert
\\
\label{sum-over-rho-new}
  &\ll
  \begin{cases}
    N & \text{if $\vert y \vert \leq 1/N$,} \\
    \vert z\vert ^{-1} + \vert z\vert ^{1/2}
    \log^2 (2 N \vert y\vert) & \text{if $\vert y \vert > 1/N$.}
  \end{cases}
\end{align}

By symmetry, we may assume that $\gamma_{1}> 0$.
By \eqref{sum-over-rho-new}, \eqref{z^-1} and \eqref{z^w}, 
for $y \in(-\infty, 0]$  we are first led to estimate
\begin{align*}
  \sum_{\rho_{1} \colon \gamma_{1} > 0} 
    \gamma_{1}^{\beta_{1} - 1/2}
    \exp\Bigl( -\frac \pi2  \gamma_{1}  \Bigr)
  &\Bigl(
    \int_{-1/N}^{0}  N^{ k +  2 + \beta_{1}} \, \dx y
    +
    \int_{-\infty}^{-1/N}  
     \frac{\dx y}{\vert y\vert^{k + 2+\beta_{1}}}    
    \\
  &+
  \int_{-\infty}^{-1/N}  \log^{2} (2N \vert y \vert)
    \frac{\dx y}{\vert y\vert^{k + 1/2+\beta_{1}}}     
  \Bigr)
  \ll_k
  N^{k + 2},
\end{align*}
by the same argument used in the proof of Lemma
\ref{series-int-zeros-alt-sign} with $\alpha=k+1/2$ and $a=1/N$.

On the other hand, for $y > 0$ we split the range of integration into
$(0, 1 / N] \cup (1 / N, +\infty)$.
By \eqref{sum-over-rho-new}, \eqref{z^-1} and Lemma
\ref{series-int-zeros-alt-sign} with $\alpha=k+1$ and $a=1/N$, on the
first interval we have
\begin{align*}
N  \sum_{\rho_{1} \colon \gamma_{1} > 0}
   \gamma_{1} ^{\beta_{1} - 1/2}
    \int_0^{1 / N} \exp\Bigl( \gamma_{1} (\arctan(N y) -\frac \pi2) \Bigr)
    \frac{\dx y}{\vert z\vert^{k + 1+\beta_{1}}}   
  &\ll_k
  N^{k + 2}.
\end{align*}
On the other interval, again by  \eqref{z^-1}, we have to estimate
\begin{align*}
  &
  \sum_{\rho_{1} \colon \gamma_{1} > 0}
   \gamma_{1} ^{\beta_{1} - 1/2}
    \int_{1 / N}^{+\infty} \exp\Bigl( - \gamma_{1} \arctan\frac 1{N y} \Bigr) 
      \frac{ y ^{-1}+ y^{1/2}\log^{2} (2Ny)}{y^{k + 1+\beta_{1}}} \dx y
 \\
  &=
  N^{k}
  \sum_{\rho_{1} \colon \gamma_{1} > 0}
  N^{\beta_{1} }
   \gamma_{1} ^{\beta_{1} - 1/2}
    \int_{1 }^{+\infty} \exp\Bigl( - \gamma_{1} \arctan\frac 1{u} \Bigr)
      \frac{Nu^{-1}+ u^{1/2}N^{-1/2}\log^{2} (2u)}{u^{k + 1+\beta_{1}}} \dx y.
\end{align*}
Lemma~\ref{series-int-zeros} with $\alpha=k+1/2$ shows that the last
term is $\ll_k N^{k + 2}$.
This implies that the integral in \eqref{conv-integral-2} is $\ll_k N^{k + 2}$ 
provided that $k > 1$ and hence we can exchange the first summation 
with the integral in this case.

To exchange the second summation we have to consider
\begin{equation}
\label{conv-integral-3}
  \sum_{\rho_{1}}
    \vert \Gamma(\rho_{1})\vert  
  \sum_{\rho_{2}}
    \vert \Gamma(\rho_{2})\vert    
      \int_{(1/N)} \vert 
    e^{N z}\vert  \vert z\vert^{- k - 1} 
    \vert z^{- \rho_{1}} \vert 
    \vert z^{- \rho_{2}} \vert \, \vert \dx z \vert.
\end{equation}
By symmetry, we can consider  $\gamma_{1},\gamma_{2}> 0$ or
$\gamma_{1} >0$, $\gamma_{2}< 0$.

Assuming $\gamma_{1},\gamma_{2}> 0$, 
for $y \leq  0$ we have
$\gamma_{j} \arctan(N y) -\frac \pi2  \gamma_{j}  \le - \frac \pi2  \gamma_{j}$,
$j=1,2$, and, by \eqref{z^w}, the corresponding contribution to
\eqref{conv-integral-3} is
\begin{align*}
  \sum_{\rho_{1} \colon \gamma_{1} > 0}
  \gamma_{1} ^{\beta_{1} - 1/2}
   & \exp\Bigl( -\frac \pi2  \gamma_{1}  \Bigr)
  \sum_{\rho_{2} \colon \gamma_{2} > 0} 
  \gamma_{2} ^{\beta_{2} - 1/2}\exp\Bigl( -\frac \pi2  \gamma_{2}  \Bigr)
  \Bigl(
    \int_{-\infty}^{0} 
     \frac{\dx y}{\vert z\vert^{k + 1+\beta_{1}+\beta_{2}}}   
    \Bigr)
    \\
  &\ll_k
  N^{k + 2 }
   \sum_{\rho_{1} \colon \gamma_{1} > 0} 
    \gamma_{1} ^{\beta_{1} - 1/2} \exp\Bigl( -\frac \pi2  \gamma_{1}  \Bigr)
   \sum_{\rho_{2} \colon \gamma_{2} > 0} 
    \gamma_{2} ^{\beta_{2} - 1/2}\exp\Bigl( -\frac \pi2  \gamma_{2}  \Bigr)
  \ll_k
  N^{k + 2},
\end{align*}
using standard zero-density estimates  and \eqref{z^-1}.
On the other hand, for $y > 0$ we split the range of integration into
$(0, 1 / N] \cup (1 / N, +\infty)$.
On the first interval we have
\begin{align*}
  \sum_{\rho_{1} \colon \gamma_{1} > 0}
   \gamma_{1} ^{\beta_{1} - 1/2}
 & \sum_{\rho_{2} \colon \gamma_{2} > 0}
  \gamma_{2} ^{\beta_{2} - 1/2}
  \int_0^{1 / N}
    \exp\Bigl( (\gamma_{1}+\gamma_{2}) (\arctan(N y) -\frac \pi2) \Bigr)
       \frac{\dx y}{\vert z\vert^{k + 1+\beta_{1}+\beta_{2}}}    \\
  &\ll
  \sum_{\rho_{1} \colon \gamma_{1} > 0}
   \gamma_{1} ^{\beta_{1} - 1/2}
  \sum_{\rho_{2} \colon \gamma_{2} > 0} 
   \gamma_{2} ^{\beta_{2} - 1/2}
   \exp\Bigl( -\frac \pi4  (\gamma_{1}+\gamma_{2})  \Bigr)
  \int_0^{1 / N}  N^{k + 1 +\beta_{1} +\beta_{2}} \, \dx y \\
  &\ll_k
   N^{k + 2}
  \sum_{\rho_{1} \colon \gamma_{1} > 0} 
  \gamma_{1} ^{\beta_{1} - 1/2} \exp\Bigl(  -\frac \pi4  \gamma_{1}  \Bigr)
  \sum_{\rho_{2} \colon \gamma_{2} > 0} 
  \gamma_{2} ^{\beta_{2} - 1/2} \exp\Bigl( -\frac \pi4  \gamma_{2}  \Bigr)
  \ll_k
  N^{k + 2},
\end{align*}
arguing as above. With similar computations,
on the other interval we have
\begin{align*}
  &
  \sum_{\rho_{1} \colon \gamma_{1} > 0}
  \gamma_{1} ^{\beta_{1} - 1/2}
  \sum_{\rho_{2} \colon \gamma_{2} > 0}
  \gamma_{2} ^{\beta_{2} - 1/2}
  \int_{1 / N}^{+\infty}
    \exp\Bigl( (\gamma_{1}+\gamma_{2}) (\arctan(N y) -\frac \pi2) \Bigr)
    \frac{\dx y}{y^{k + 1+\beta_{1}+\beta_{2}}} \\
  &=
  N^{k}
  \sum_{\rho_{1} \colon \gamma_{1} > 0}
  N^{\beta_{1}}
  \gamma_{1} ^{\beta_{1} - 1/2}
  \sum_{\rho_{2} \colon \gamma_{2} > 0}
  N^{\beta_{2}}
  \gamma_{2} ^{\beta_{2} - 1/2}
    \int_{1}^{+\infty} \exp\Bigl(- (\gamma_{1}+\gamma_{2}) \arctan\frac 1u \Bigr)
      \frac{\dx u}{u^{k + 1+\beta_{1}+\beta_{2}}}.     
\end{align*}
Arguing as in the proof of Lemma~\ref{series-int-zeros}, the integral
on the right is $\sameorder (\gamma_1 + \gamma_2)^{-k-\beta_1-\beta_2}$.
The inequality
\begin{equation}
\label{ineq-zeros}
  \frac{\gamma_1^{\beta_1 - 1/2} \gamma_2^{\beta_2 - 1/2}}
       {(\gamma_1 + \gamma_2)^{\beta_1 + \beta_2}}
  \le
  \frac1{\gamma_1^{1/2} \gamma_2^{1/2}}
\end{equation}
shows that it is sufficient to consider
\begin{align*}
  N^{k}
  \sum_{\rho_{1} \colon \gamma_{1} > 0}
  \sum_{\rho_{2} \colon \gamma_{2} > 0}
    N^{\beta_{1}+\beta_{2}}
    \frac1{\gamma_{1} ^{1/2}\gamma_{2} ^{1/2}(\gamma_1 + \gamma_2)^{k}}  
  &\ll
  N^{k+2}
  \sum_{\rho_{1} \colon \gamma_{1} > 0}
    \frac1{\gamma_1^{k + 1/2}}
    \sum_{\rho_{2} \colon  0< \gamma_2 \le \gamma_1}
      \frac1{\gamma_2^{1/2}}
  \\
  &\ll
  N^{k+2}
  \sum_{\rho_{1} \colon \gamma_{1} > 0}
   \frac{\log \gamma_1}{\gamma_1^{k}}
\end{align*}
and the last series over zeros converges for $k > 1$.
 
Assume now $\gamma_{1}>0$, $\gamma_{2}< 0$.
For $y \leq  0$ we have
$\gamma_{1} \arctan(N y) -\frac \pi2  \gamma_{1}  \le - \frac \pi2  \gamma_{1}$,  
by \eqref{z^-1} the corresponding contribution to \eqref{conv-integral-3} is
\begin{align*} 
\ll_{k}    
 &
  \sum_{\rho_{1} \colon \gamma_{1} > 0}  \gamma_{1} ^{\beta_{1} - 1/2} 
   \exp\Bigl( -\frac \pi2  \gamma_{1}  \Bigr)
  \Bigl\{
   \sum_{\rho_{2} \colon \gamma_{2} < 0}   
     \vert \gamma_{2} \vert ^{\beta_{2} - 1/2} 
     \Bigl[
     \exp\Bigl( -\frac \pi4 \vert  \gamma_{2}\vert   \Bigr)
    \int_{-1/N}^{0}
    N^{k + 1 +\beta_{1}+\beta_{2}} \, \dx y 
    \\
    &
    \hskip 0.5cm
    +    
    \int_{-\infty}^{-1/N} 
    \exp\Bigl(- \vert  \gamma_{2}\vert ( \arctan(N y) +\frac \pi2) \Bigr)
    \frac{\dx y}{\vert y\vert ^{k + 1+\beta_{1}+\beta_{2}}}  
     \Bigr]
  \Bigr\}  \\
  &\ll_k
  N^{k+2} 
   \sum_{\rho_{1} \colon \gamma_{1} > 0}  
   \gamma_{1} ^{\beta_{1} - 1/2}  
   \exp\Bigl( -\frac \pi2  \gamma_{1}  \Bigr)
\sum_{\rho_{2} \colon \gamma_{2} < 0} 
 \vert \gamma_{2} \vert ^{\beta_{2} - 1/2}  
 \exp\Bigl( -\frac \pi4 \vert \gamma_{2} \vert  \Bigr)
\\
&
+ 
 N^{k + 2 }
 \!\! \sum_{\rho_{1} \colon \gamma_{1} > 0} \!\! \gamma_{1} ^{\beta_{1} - 1/2}  
  \exp\Bigl( -\frac \pi2  \gamma_{1}  \Bigr)
 \!\! \sum_{\rho_{2} \colon \gamma_{2} < 0} \!\!
  \vert \gamma_{2} \vert ^{\beta_{2} - 1/2}  
  \int_1^{+\infty} \!\!\!\!
   \exp\Bigl( - \vert \gamma_{2}\vert \arctan\frac 1u \Bigr)
      \frac{\dx u}{u^{k + 1 +\beta_{1}+\beta_{2}}} 
    \\
  &\ll_k
  N^{k + 2}
  +
  N^{k + 2}
  \sum_{\rho_{1} \colon \gamma_{1} > 0} 
  \gamma_{1} ^{\beta_{1} - 1/2}  
  \exp\Bigl( -\frac \pi2 \gamma_{1} \Bigr) 
  \ll_k
  N^{k + 2}
\end{align*}
for $k > 1/2$, by Lemma~\ref{series-int-zeros} and standard zero-density
estimates.

On the other hand, the case $\gamma_{1}>0$, $\gamma_{2}< 0$ and $y > 0$
can be estimated in a similar way essentially exchanging the role of
$\gamma_{1}$ and $\gamma_{2}$ in the previous argument.

This implies that the integral in \eqref{conv-integral-3} is $\ll_k N^{k + 2}$ 
provided that $k > 1$. Combining the convergence conditions for 
\eqref{conv-integral-2}-\eqref{conv-integral-3},
we see that we can exchange both summations with the integral provided that $k>1$.

\section{Convergence of the double sum over zeros}
\label{sec:double-sum}

In this section we prove that the double sum on the right of
\eqref{expl-form-Goldbach-bis} converges absolutely 
for every $k > 1 / 2$.
We need \eqref{Stirling} uniformly for $x \in [0, k + 3]$ and
$\vert y \vert \ge T$, where $T$ is large but fixed: this provides
both an upper and a lower bound for $\vert \Gamma(x + i y) \vert$.
Let
\begin{equation*}
  \Sigma 
  =
  \sum_{\rho_1} \sum_{\rho_2}
    \Bigl\vert \frac{\Gamma(\rho_1) \Gamma(\rho_2)}
                {\Gamma(\rho_1 + \rho_2 + k + 1)}
    \Bigr\vert,
\end{equation*}
so that, by the symmetry of the zeros of the Riemann zeta-function, we
have
\begin{align*}
  \Sigma
  &=
  2
  \sum_{\rho_{1} \colon \gamma_{1} > 0} \sum_{\rho_{2} \colon \gamma_2 > 0}
    \Bigl\vert
      \frac{\Gamma(\rho_1) \Gamma(\rho_2)}
           {\Gamma(\rho_1 + \rho_2 + k + 1)}
    \Bigr\vert
  +
  2
  \sum_{\rho_{1} \colon \gamma_{1} > 0} \sum_{\rho_{2} \colon \gamma_2 > 0}
    \Bigl\vert
      \frac{\Gamma(\rho_1) \Gamma(\overline{\rho}_2)}
           {\Gamma(\rho_1 + \overline{\rho}_2 + k + 1)}
    \Bigr\vert \\
  &=
  2 (\Sigma_1 + \Sigma_2),
\end{align*}
say.
It is clear that if both $\Sigma_1$ and $\Sigma_2$ converge, then the
double sum on the right-hand side of \eqref{expl-form-Goldbach-bis}
converges absolutely.
In order to estimate $\Sigma_1$ we choose a large $T$ and let
\begin{align*}
  D_0
  &=
  \{ (\rho_1, \rho_2)\colon (\gamma_1, \gamma_2) \in [0, 2 T]^2 \}, \\
  D_1
  &=
  \{ (\rho_1, \rho_2)\colon 
    \gamma_1 \ge T, \,
    T \le \gamma_2 \le \gamma_1 \}, \\
  D_2
  &=
  \{ (\rho_1, \rho_2)\colon 
    \gamma_1 \ge T, \,
    0 \le \gamma_2 \le T \}, \\
  D_3
  &=
  \{ (\rho_1, \rho_2)\colon 
    \gamma_2 \ge T, \,
    T \le \gamma_1 \le \gamma_2  \}, \\
  D_4
  &=
  \{ (\rho_1, \rho_2)\colon 
    \gamma_2 \ge T, \,
    0 \le \gamma_1 \le T  \},
\end{align*}
so that
$\Sigma_1 \le \Sigma_{1,0} + \Sigma_{1,1} + \Sigma_{1,2} + \Sigma_{1,3} +
 \Sigma_{1,4}$,
say, where $\Sigma_{1,j}$ is the sum with $(\rho_1, \rho_2) \in D_j$.
Now, $D_0$ contributes a bounded amount, that depends only on $T$,
and, by symmetry again, $\Sigma_{1,1} = \Sigma_{1,3}$ and
$\Sigma_{1,2} = \Sigma_{1,4}$.
We also recall the inequality \eqref{ineq-zeros}
which is valid for all couples of zeros considered in $\Sigma_1$.
Hence
\begin{align*}
  \Sigma_{1,1}
  &\ll
  \doublesum_{\substack{\rho_{1} \colon \gamma_1 \ge T \\ \rho_{2} \colon  T \le \gamma_2 \le \gamma_1}}
    \Bigl\vert
      \frac{\Gamma(\beta_1 + i \gamma_1) \Gamma(\beta_2 + i \gamma_2)}
           {\Gamma(\beta_1 + \beta_2 + k + 1 + i (\gamma_1 + \gamma_2))}
    \Bigr\vert \\
  &\ll
  \doublesum_{\substack{\rho_{1} \colon \gamma_1 \ge T \\ \rho_{2} \colon T \le \gamma_2 \le \gamma_1}}
    \frac{e^{- \pi (\gamma_1 + \gamma_2) / 2} \gamma_1^{\beta_1 - 1/2} \gamma_2^{\beta_2 - 1/2}}
         {e^{- \pi (\gamma_1 + \gamma_2) / 2} (\gamma_1 + \gamma_2)^{\beta_1 + \beta_2 + k + 1/2}} 
  \ll
  \doublesum_{\substack{\rho_{1} \colon \gamma_1 \ge T \\ \rho_{2} \colon T \le \gamma_2 \le \gamma_1}}
    \frac1{\gamma_1^{1/2} \gamma_2^{1/2} (\gamma_1 + \gamma_2)^{k + 1/2}}
         \\
   & 
   \ll
  \sum_{\rho_{1} \colon \gamma_1 \ge T}
    \frac1{\gamma_1^{k + 1}}
    \sum_{\rho_{2} \colon T \le \gamma_2 \le \gamma_1}
      \frac1{\gamma_2^{1/2}}
  \ll
  \sum_{\rho_{1} \colon \gamma_1 \ge T}
   \frac{\log \gamma_1}{\gamma_1^{k + 1 / 2}}.
\end{align*}
A similar argument proves that
\[
  \Sigma_{1,2}
  \ll_{k,T}
  \sum_{\rho_{1} \colon \gamma_1 \ge T}
  \frac1{\gamma_1^{k + 1}},
\]
since $\Gamma(\rho_2)$ is uniformly bounded, in terms of $T$, for
$(\rho_1, \rho_2) \in D_2$.
Summing up, we have
\[
  \Sigma_1
  \ll_{k,T}
  1
  +
  \sum_{\rho_{1} \colon \gamma_1 \ge T}
  \frac{\log \gamma_1}{\gamma_1^{k + 1/2}},
\]
which is convergent provided that $k > 1 / 2$.
In order to estimate $\Sigma_2$ we use a similar argument.
Choose a large $T$ and let
\begin{align*}
  E_0
  &=
  \{ (\rho_1, \rho_2)\colon (\gamma_1, \gamma_2) \in [0, 2 T]^2 \}, \\
  E_1
  &=
  \{ (\rho_1, \rho_2)\colon 
   \gamma_1 \ge 2 T, \,
    0 \le \gamma_2 \le T \}, \\
  E_2
  &=
  \{ (\rho_1, \rho_2)\colon 
    \gamma_1 \ge 2 T, \,
    T \le \gamma_2 \le \gamma_1 - T \}, \\
  E_3
  &=
  \{ (\rho_1, \rho_2)\colon 
    \gamma_1 \ge 2 T, \,
    \gamma_1 - T \le \gamma_2 \le \gamma_1 \}, \\
  E_4
  &=
  \{ (\rho_1, \rho_2)\colon  
    \gamma_2 \ge 2 T, \,
    \gamma_2 - T \le \gamma_1 \le \gamma_2 \}, \\
  E_5
  &=
  \{ (\rho_1, \rho_2)\colon 
    \gamma_2 \ge 2 T, \,
    T \le \gamma_1 \le \gamma_2 - T \}, \\
  E_6
  &=
  \{ (\rho_1, \rho_2)\colon 
    \gamma_2 \ge 2 T, \,
    0 \le \gamma_1 \le T \},
\end{align*}
so that
$\Sigma_2 \le \Sigma_{2,0} +\Sigma_{2,1} + \Sigma_{2,2} + \Sigma_{2,3} +
 \Sigma_{2,4} + \Sigma_{2,5} + \Sigma_{2,6}$,
say, where $\Sigma_{2,j}$ is the sum with $(\rho_1, \rho_2) \in E_j$.
Now, $E_0$ contributes a bounded amount, that depends only on $T$,
and, by symmetry again, $\Sigma_{2,1} = \Sigma_{2,6}$,
$\Sigma_{2,2} = \Sigma_{2,5}$ and $\Sigma_{2,3} = \Sigma_{2,4}$.
Again we use \eqref{Stirling} as above; hence
\begin{align*}
  \Sigma_{2,2}
  &=
  \doublesum_{\substack{\rho_{1} \colon \gamma_1 \ge 2 T \\ \rho_{2} \colon  T \le \gamma_2 \le \gamma_1 - T}}
    \Bigl\vert
      \frac{\Gamma(\beta_1 + i \gamma_1) \Gamma(\beta_2 - i \gamma_2)}
           {\Gamma(\beta_1 + \beta_2 + 1 + i (\gamma_1 - \gamma_2))}
    \Bigr\vert 
  \ll
  \sum_{\substack{\rho_{1} \colon \gamma_1 \ge 2 T \\ \rho_{2} \colon T \le \gamma_2 \le \gamma_1 - T}}
    \frac{\gamma_1^{\beta_1 - 1/2} \gamma_2^{\beta_2 - 1/2} e^{-\pi \gamma_2}}
         {(\gamma_1 - \gamma_2)^{\beta_1 + \beta_2 + k + 1/2}} \\
  &\ll
  \sum_{\rho_{1} \colon \gamma_1 \ge 2 T}
    \gamma_1^{\beta_1 - 1/2} \log \gamma_1
    \int_T^{\gamma_1 - T}
      \frac{t^{1 / 2}}{(\gamma_1 - t)^{\beta_1 + k + 1/2}} \, e^{-\pi t} \, \dx t \\
  &\ll
  \sum_{\rho_{1} \colon \gamma_1 \ge 2 T}
    e^{-\pi \gamma_1}
    \gamma_1^{\beta_1} \log \gamma_1
    \int_T^{\gamma_1 - T}
      \frac{e^{\pi u} \, \dx u}{u^{\beta_1 + k + 1/2}} \\
  &\ll
  \sum_{\rho_{1} \colon \gamma_1 \ge 2 T}
    e^{-\pi \gamma_1}
    \gamma_1^{\beta_1} \log \gamma_1
    \frac{e^{\pi (\gamma_1 - T)}}{(\gamma_1 - T)^{\beta_1 + k + 1/2}} 
  \ll_T
  \sum_{\rho_{1} \colon \gamma_1 \ge 2 T}
    \frac{\log \gamma_1}{\gamma_1^{k + 1/2}}.
\end{align*}
The rightmost series over zeros plainly converges for $k > 1 / 2$.
The contribution of zeros in $E_1$ is treated in a similar fashion,
using the uniform upper bound $\Gamma(\rho_2) \ll_T 1$, and is
smaller.
We now deal with $\Sigma_{2,3}$: we have
\begin{align*}
  \Sigma_{2,3}
  &=
  \doublesum_{\substack{\rho_{1} \colon \gamma_1 \ge 2 T \\ \rho_{2} \colon  \gamma_1 - T \le  \gamma_2 \le \gamma_1}}
    \Bigl\vert
      \frac{\Gamma(\beta_1 + i \gamma_1) \Gamma(\beta_2 - i \gamma_2)}
           {\Gamma(\beta_1 + \beta_2 + k + 1 + i (\gamma_1 - \gamma_2))}
    \Bigr\vert \\
  &\ll
  \sum_{\rho_{1} \colon \gamma_1 \ge 2 T}
    e^{- \pi \gamma_1 / 2} \gamma_1^{\beta_1 - 1/2}
    \sum_{\rho_{2} \colon \gamma_1 - T \le  \gamma_2 \le \gamma_1}
      e^{- \pi \gamma_2 / 2} \gamma_2^{\beta_2 - 1/2}
      \Bigl(
        \min_{\substack{k + 1 \le x \le k + 3 \\ 0 \le t \le T}} \vert \Gamma(x + i t) \vert
      \Bigr)^{-1} \\
  &\ll_{k, T}
  \sum_{\rho_{1} \colon \gamma_1 \ge 2 T}
    e^{- \pi \gamma_1} \gamma_1^{\beta_1 + 1}
    \log(\gamma_1 + T),
\end{align*}
provided that $T$ is large enough.
Here we are using Theorem 9.2 of Titchmarsh \cite{Titchmarsh1986} with
$T$ large but fixed.
The series at the extreme right is plainly convergent.


\providecommand{\bysame}{\leavevmode\hbox to3em{\hrulefill}\thinspace}
\providecommand{\MR}{\relax\ifhmode\unskip\space\fi MR }
\providecommand{\MRhref}[2]{%
  \href{http://www.ams.org/mathscinet-getitem?mr=#1}{#2}
}
\providecommand{\href}[2]{#2}

\vskip0.5cm
\noindent
\begin{tabular}{l@{\hskip 24mm}l}
Alessandro Languasco               & Alessandro Zaccagnini\\
Universit\`a di Padova     & Universit\`a di Parma\\
Dipartimento di Matematica & Dipartimento di Matematica \\
Via Trieste 63                & Parco Area delle Scienze, 53/a \\
35121 Padova, Italy            & 43124 Parma, Italy\\
{\it e-mail}: languasco@math.unipd.it        & {\it e-mail}:
alessandro.zaccagnini@unipr.it  
\end{tabular}

\end{document}